# Optimal Local Multi-scale Basis Functions for Linear Elliptic Equations with Rough Coefficient*

Thomas Y. Hou        Pengfei Liu

**Abstract**

This paper addresses a multi-scale finite element method for second order linear elliptic equations with arbitrarily rough coefficient. We propose a local oversampling method to construct basis functions that have optimal local approximation property. Our methodology is based on the compactness of the solution operator restricted on local regions of the spatial domain, and does not depend on any scale-separation or periodicity assumption of the coefficient. We focus on a special type of basis functions that are harmonic on each element and have optimal approximation property. We first reduce our problem to approximating the trace of the solution space on each edge of the underlying mesh, and then achieve this goal through the singular value decomposition of an oversampling operator. Rigorous error estimates can be obtained through thresholding in constructing the basis functions. Numerical results for several problems with multiple spatial scales and high contrast inclusions are presented to demonstrate the compactness of the local solution space and the capacity of our method in identifying and exploiting this compact structure to achieve computational savings.

**Keywords.** Multi-scale finite element method. Oversampling. Optimal local basis. High-contrast.

## 1 Introduction

Many problems of practical importance in science and engineering have multi-scale feature: composite materials modeling and flows in porous media are typical examples of such kind. In many cases, the qualities of interest are only related to the large-scale properties of the solutions. However, the fine-scale features of the model can have significant impact on the large-scale properties of the solutions, thus one needs to use very fine mesh to resolve the small-scale variations of the problem even though only large-scale properties of the solutions are required. The computational cost can be prohibitive. For these so-called multi-scale problems, efficient upscaling methods that allow us to incorporate the small-scale features of the problem into the large-scale properties of the solutions are desired.

In this work, we use the following example of second order linear elliptic equation with homogeneous Dirichlet boundary condition to illustrate our upscaling methodology,

$$\begin{cases} -\text{div}(a(x)\nabla u(x)) = f(x), & x \in \Omega; \\ u(x)|_{\partial\Omega} = 0, \end{cases} \quad (1.1)$$

where $\Omega$ is a convex polygon domain in $R^d$ with $d = 2, 3$. We assume that the equation is uniformly elliptic, i.e., there exist $\lambda_{\min} > 0$ and $\lambda_{\min} > 0$ such that

$$a(x) \in [\lambda_{\min}, \lambda_{\max}]. \quad (1.2)$$

We do not assume any regularity of the coefficient $a(x) \in L^{\infty}(\Omega)$, which may have multiple spatial scales, thus the above equation (1.1) can be used to model diffusion process in strongly heterogeneous media. We also assume that in (1.1) the forcing function $f(x) \in L^2(\Omega)$. Then the existence of solution to (1.1), $u(x) \in H_0^1(\Omega)$ follows immediately from the Lax-Milgram theorem, and we have

$$c\|f\|_{H^{-1}(\Omega)} \leq \|u(x)\|_{H_0^1(\Omega)} \leq C\|f\|_{H^{-1}(\Omega)}. \quad (1.3)$$

However, due to the roughness of the coefficient $a(x)$, the solution to (1.1) $u(x)$ loses regularity and ceases to be in $H^2(\Omega)$. Classical finite element method uses piecewise polynomials to approximate the solution space, and its convergence depends on the following approximation property

$$\|u(x) - Ju(x)\|_{H_0^1(\Omega)} \leq Ch\|u(x)\|_{H^2(\Omega)}, \quad (1.4)$$

---

and the regularity result
$$\|u(x)\|_{H^2(\Omega)} \leq C\|f(x)\|_{L^2(\Omega)}, \tag{1.5}$$

where $Ju$ is the piecewise polynomial interpolation of $u(x)$, and $h$ is the underlying mesh size. Classical finite element method may fail for these problems (1.1), since $\|u(x)\|_{H^2}$ cannot be bounded by $\|f(x)\|_{L^2(\Omega)}$ in (1.5). It is actually showed in [6] that the polynomial finite elements can perform arbitrarily badly in this setting, thus (1.1) can serve as a typical example of multi-scale problems.

One of the strategies to numerically solve the multi-scale problem (1.1) when classical finite element methods fail is using problem-dependent basis that incorporates properties of the coefficient $a(x)$, to approximate the solution space. To be specific, one constructs basis functions

$$\phi_1(x), \phi_2(x), \ldots \phi_n(x) \in H_0^1(\Omega), \tag{1.6}$$

that may depend on the elliptic operator $-\text{div}(a(x)\nabla(\cdot))$, and find numerical solution

$$u_h(x) \in V_h(x) = \text{span}\{\phi_1(x), \phi_2(x), \ldots \phi_n(x)\} \subset H_0^1(\Omega), \tag{1.7}$$

using the Galerkin projection. Namely, we find $u_h(x)$ within the trial space $V_h$, such that

$$a(u_h(x), v(x)) = \langle f(x), v(x) \rangle, \quad \text{for all} \quad v \in V_h, \tag{1.8}$$

where

$$a(u(x), v(x)) = \int_\Omega \nabla u(x)^t a(x) \nabla v(x) \mathrm{d}x, \quad \langle f(x), v(x) \rangle = \int_\Omega f(x) v(x) \mathrm{d}x. \tag{1.9}$$

The numerical solution defined above satisfies the following optimal property under the energy norm

$$\|u(x) - u_h(x)\|_a = \inf_{v(x) \in V_h} \|u(x) - v(x)\|_a, \tag{1.10}$$

where the energy norm is equivalent to the $H_0^1(\Omega)$ norm, and defined as

$$\|u(x)\|_a^2 = a(u(x), u(x)) = \int_\Omega \nabla u(x)^t a(x) \nabla u(x) \mathrm{d}x. \tag{1.11}$$

In this work we will employ the above strategy to numerically solve (1.1). Note that to obtain the numerical solution $u_h(x)$ from the Galerkin projection (1.8), one needs to solve a linear system of size $n \times n$. Thus to make the computational cost small, we want the number of the basis functions used in (1.6) to be small. Besides, we want the basis functions in (1.6) to have compact support such that the stiffness matrix formed in (1.8) is sparse thus easy to compute and invert.

We propose an effective method to construct basis functions (1.6) with optimal local approximation property in this paper. Our method is based on the compactness of the solution operator to (1.1) restricted on local regions of the domain. To be specific, we introduce the following operator

$$T_i: \quad f(x) \to u_i(x) = u(x)|_{D_i}, \tag{1.12}$$

where $D_i$ is local subset of $\Omega$ of size $O(H)$, and $H$ is chosen according to the desired order of accuracy. We construct local basis functions that can approximate the range of $T_i$ with controlled accuracy, and combine them together to get the approximate solution space to (1.1). The compactness of the operator $T_i$, after removing the singular part, will be demonstrated numerically in section 2.

On each local region of the domain, $D_i$, we divide the local solution $u_i(x)$ into two orthogonal parts with respect to energy norm (1.9): an $a(x)$-harmonic part, and a local bubble part. We show that the bubble part of the solution is small and its compact structure can be easily identified by inverting the elliptic equation (1.1) locally on each region $D_i$. We consider approximating the $a(x)$-harmonic part of the solution space using a special type of basis functions that are $a(x)$-harmonic on each $D_i$ (but not across the boundary of $D_i$), and call basis functions of such type **multi-scale basis**. Due to the smallness of the bubble part of the solution, we demonstrate that multi-scale basis functions are optimal in approximating the solution space for fixed local boundary conditions on $\partial D_i$.

The $a(x)$-harmonic part of the solution only depends on the restriction of the solution on the boundary of the local regions $D_i$, and we seek to identify the compact structure of the trace of the solution space on $\partial D_i$. Using a primary set of interpolation basis functions, $\psi_i(x)$, we can reduce our problem to approximating the solution space on each edge of the coarse mesh, $e$. We introduce an oversampling operator that maps the solution on an oversampling domain $W$ to the solution restricted on an edge of the mesh $e$. Then we employ the compactness of the oversampling operator to construct the optimal boundary basis functions on each edge $e$. The optimal choice of the interpolation basis functions $\psi_i(x)$ can also be identified by solving least square problems. With these local boundary basis functions, we construct basis functions (1.6) that approximate the $a(x)$-harmonic part of the solution space by solving some local boundary value problems. The resulting basis functions (1.6) are $a(x)$-orthogonal with respect to the bubble part of the solution space, and because of this property we can add the bubble part back to our numerical solutions by simply solving some local cell problems.

The resulting method consists of two stages: in the offline stage we identify the local compact structure of the solution space, and build multi-scale basis functions and the corresponding stiffness matrix in (1.8); in the online stage, for any given forcing function $f(x) \in L^2(\Omega)$, we solve the equation (1.1) efficiently using the multi-scale basis functions constructed offline with a very low computation cost. Our method can achieve significant computational savings in the multi-query setting where equation (1.1) need to be solved for multiple times with different forcing functions $f(x) \in L^2(\Omega)$.

Numerical examples for several problems with rough coefficients and high-contrast channels are presented. Our method achieves high accuracy for these problems, and these numerical results suggest that our method is very robust and works equally well for problems without scale separation, or have high contrast channels. We demonstrate our methodology through the second-order scalar elliptic equation, but it can be applied to other linear elliptic problems like elasticity equations without much difficulty.

There is a fast growing literature on multi-scale methods for elliptic equations with rough coefficient, and we briefly review some related work below. The classical homogeneziation theory [7] considers $a(x)$ with periodic structure $a(x) = a(x, \frac{x}{\epsilon})$, and derives an effective equation governing the asymptotic behaviors of the solutions as the small scale $\epsilon$ goes to 0. In [21, 18, 22, 14], the authors proposed the Multi-scale finite element method (MsFEM), in which the multi-scale basis functions are constructed by solving local elliptic problems with homogeneous right hand side, and the convergence analysis of MsFEM in the periodic setting was given in [18, 15]. An oversampling technique [18] was proposed to reduce the resonance error introduced due to the artificial boundary conditions in solving the local problems. This present work can be viewed as a continuation of the MsFEM, and we propose a robust method to choose and enrich the local boundary conditions using the oversampling operator. In [27], the authors showed that the solutions gain an order of regularity with respect to the Harmonic coordinate [3], and constructed multi-scale basis functions using this property. The Harmonic coordinate was recently employed in [11] for global upscaling. In [8, 28], the flux norm was introduced and employed to show the compactness of the solution space and construct (localized) basis functions. In [29] the polyharmonic spline was employed, which was later put in the Bayesian inference setting [26]. In [24], the generalized finite element method was proposed, which provides a general framework to combine local approximation spaces together using a partition of unity formulation. In [5, 4], the local special basis functions are constructed by solving some local spectral problems, and then combined together using the partition of unity framework. The generalized multi-scale finite element method [12] is a systematic method to construct multi-scale basis functions for a family multi-scale problems with parameters. The Heterogeneous multi-scale method [1, 30] numerically decomposes the structure of the medium into a micro-scale and a macro-scale, and compute solutions of cell problems on the micro-scale to get the local homogenized matrices. Finally we remake that multi-scale methods using problem-dependent basis functions were also employed to solve problems in which the coefficient has high contrast inclusions [9].

The remaining part of this paper is organized as follows. In section 2, we demonstrate the compactness of the solution operator restricted on local regions of the spatial domain. In section 3, we divide the solutions on each local region of the domain to different parts corresponding to the trace of the solution on the edges of the coarse mesh, and identify their compact structures separately. And with this compact structure, we construct multi-scale basis functions that have optimal local approximation property. In section 4, numerical results are presented to demonstrate the capacity of our method in identifying and exploiting the compactness of the solution space to achieve computational savings. Concluding remarks are made in section 5.

## 2 Compactness of the Solution Space Restricted on Local Regions of the Domain

The existence of finite number of basis functions (1.6) that can approximate the solutions space to (1.1) up to any accuracy is implied by the compactness of the solution operator, $T$, which maps from the forcing function $f(x) \in L^2(\Omega)$ to the corresponding solution $u(x) \in H_0^1(\Omega)$.

$$T: \quad f(x) \in L^2(\Omega) \to u(x) \in H_0^1(\Omega). \tag{2.1}$$

The compactness of $T$ is demonstrated in [23, 8], and it was employed for elliptic equations with random input data recently in [19, 20] for stochastic model reduction.

To be specific, the solution operator $T$ can be decomposed as

$$T = L^{-1} I_{L^2(\Omega) \to H^{-1}(\Omega)}, \tag{2.2}$$

where $L^{-1}$ maps $f(x) \in H^{-1}(\Omega)$ to the solution $u(x) \in H_0^1(\Omega)$, and $I_{L^2(\Omega) \to H^{-1}(\Omega)}$ is the embedding operator from $L^2(\Omega)$ to $H^{-1}(\Omega)$. From (1.3), we can see that $L^{-1}$ is continuous and indeed a homomorphism, and the compactness of $I_{L^2(\Omega) \to H^{-1}(\Omega)}$ is well known based on the Sobolev space theory [16]. Thus the compactness of $T$ follows from the decomposition (2.2). To quantify the approximability of $T$ by a finite-rank operator, we consider its Kolmogorov-$n$ width, which is defined in below.

**Definition 2.1** (Kolmogorov $n$-width). For a compact linear operator $T$ that maps between two Hilbert spaces, we define its Kolmogorov $n$-width as

$$d_n(T) = \inf_{T_n} \|T - T_n\|, \tag{2.3}$$

where $T_n$ runs over all rank-$n$ linear operators.

Due to the fact that $L^{-1}$ is a homomorphism, one can easily see that the Komogorov-$n$ width of $T$ is only different from that of $I_{L^2(\Omega) \to H^{-1}(\Omega)}$ by a factor that depends on $\lambda_{\min}$, $\lambda_{\max}$ (1.2) and $\Omega$,

$$c d_n(I_{L^2(\Omega) \to H^{-1}(\Omega)}) \leq d_n(T) \leq C d_n(I_{L^2(\Omega) \to H^{-1}(\Omega)}). \tag{2.4}$$

The Kolmogorov-$n$ width of the embedding operator is well-known [23, 8], and we have

$$d_n(I_{L^2(\Omega) \to H_0^1(\Omega)}) = n^{-1/d}(C + o(1)), \quad n \to \infty. \tag{2.5}$$

From (2.4) (2.5) and (2.3), we obtain that there exist of a set of basis functions, (1.6), with the following approximation property to the solution space of (1.1),

$$\sup_{\|f(x)\|_{L^2(\Omega)}=1} \inf_{c_i} \|\sum_{i=1}^n c_i \phi_i(x) - u(x)\|_{H_0^1(\Omega)} \leq C n^{-1/d}. \tag{2.6}$$

The approximation property (2.6) is optimal, and does not depend on the regularity of the coefficient $a(x)$. For practical applications in multi-scale problems, we want the basis functions $\phi_i(x)$ to have local support such that the corresponding stiffness matrix in (1.8) is sparse and easy to invert. However, the basis functions in (2.6) whose existence is implied by (2.4), (2.5) and (2.6) may be nonlocal.

Since our objective is finding basis functions (1.6) with local support, we consider a local region of the domain, $D$ with diameter $O(H)$, and a slightly larger local domain which contains $D$, $W$, which we call the oversampling region. We consider the restriction of the solutions to (1.1) on $W$,

$$u_W(x) = u(x)|_W. \tag{2.7}$$

The local solution $u_W(x)$ can be divided into two parts,

$$u_W(x) = u_W^1(x) + u_W^2(x), \tag{2.8}$$

where

$$\begin{cases} -\mathrm{div}(a(x)\nabla u_W^1(x)) = 0, & x \in W, \\ u_W^1(x) = u_W(x), & x \in \partial W, \end{cases} \tag{2.9}$$

and

$$\begin{cases} -\mathrm{div}(a(x)\nabla u_W^2(x)) = f(x), & x \in W, \\ u_W^2(x) = 0, & x \in \partial W. \end{cases} \tag{2.10}$$

We call the first part $u_W^1(x)$ the local $a(x)$-harmonic part, and the second part $u_W^2(x)$ the local bubble part. The two parts are orthogonal with respect the local inner product, $a_W(\cdot, \cdot)$,

$$a_W(u_W^1(x), u_W^2(x)) = \int_W \nabla u_W^1(x)^t a(x) \nabla u_W^2(x) \mathrm{d}x = 0. \tag{2.11}$$

The local bubble part $u_W^2(x)$ is small in the sense that

$$\|u_W^2(x)\|_{H_0^1(W)}^2 \leq C H^2 \|f(x)\|_{L^2(W)}^2, \tag{2.12}$$

which can be obtained from (1.3) and a scaling argument. (2.12) implies that if we only want to obtain $O(H)$ accuracy in our numerical solution, we can simply neglect the local bubble part.

Then we consider a local solution operator $T_D$ that maps $f(x)$ to the local $a(x)$-harmonic part, $u_W^1(x)$ restricted on $D$,

$$T_D: \ f(x) \in L^2(\Omega) \to u_W^1(x)|_D \in H^1(D), \tag{2.13}$$

and we want to construct local basis functions on $D$ that can approximate the range of $T_D$. To demonstrate the compactness of $T_D$, we choose a set of orthonormal basis in the domain and range of $T_D$ to discretize $T_D$ as a matrix, and compute the decay of its singular values. We consider the following choice of coefficient in (1.1), which has multiple fine spatial scales and is illustrated in Figure 1a,

$$a(x) = \frac{1}{6}\left(\frac{1.1 + \sin(2\pi x/\epsilon_1)}{1.1 + \sin(2\pi y/\epsilon_1)} + \frac{1.1 + \sin(2\pi y/\epsilon_2)}{1.1 + \cos(2\pi x/\epsilon_2)} + \frac{1.1 + \cos(2\pi x/\epsilon_3)}{1.1 + \sin(2\pi y/\epsilon_3)} \right.$$
$$\left. + \frac{1.1 + \sin(2\pi y/\epsilon_4)}{1.1 + \cos(2\pi x/\epsilon_4)} + \frac{1.1 + \cos(2\pi x/\epsilon_5)}{1.1 + \sin(2\pi y/\epsilon_5)} + \sin(4x^2 y^2) + 1 \right), \tag{2.14}$$

where $\epsilon_1 = \frac{1}{5}$, $\epsilon_2 = \frac{1}{13}$, $\epsilon_3 = \frac{1}{17}$, $\epsilon_4 = \frac{1}{31}$, $\epsilon_5 = \frac{1}{65}$.

We choose $\Omega = [0,1] \times [0,1]$, the oversampling region $W = [14H, 17H] \times [14H, 17H]$, and the local region $D = [15H, 16H] \times [15H, 16H]$, where $H = 1/32$. The decay of the singular values of the local solution operator (2.13) is plotted in Figure 1b. Then we compute the singular values for the local solution operator (2.13) to the Possion equation in the same setting, the decay of which is plotted in Figure 1c. From the figures, we can see that the singular values of the local solution operator decay very fast, and this fast decay does not deteriorate due to the roughness of the coefficient.

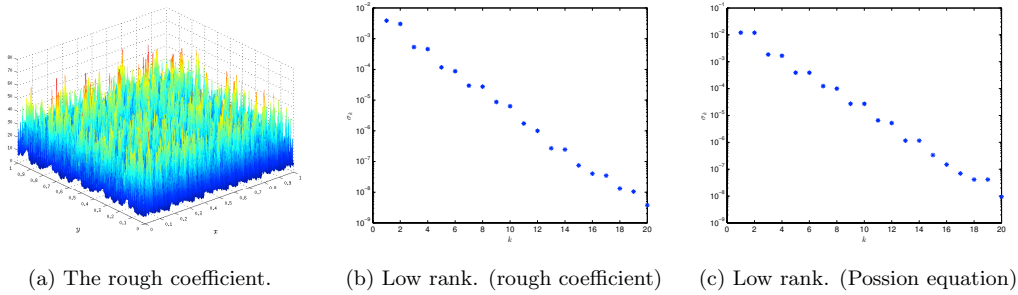

(a) The rough coefficient.  (b) Low rank. (rough coefficient)  (c) Low rank. (Possion equation)

Figure 1

The fast decay of singular values of $T_D$ implies that we can use a very small number of local basis functions, to be specific, the first several left singular vectors of $T_D$, to get very good local approximation property. However, we cannot afford to construct $T_D$ explicitly since it is a solution operator and its construction involves solving the equation (1.1) a large amount of times (globally). It is known that for a low-rank operator, the main action of $T_D$ can be captured in its image on some random vectors. This, to some degree, explains the success of some global upscaling methods [13, 11, 27] that use sampled global solutions to (1.1) to approximate the solution space to (1.1) locally.

We will not pursue this perspective in this work. Instead, we decompose $T_D$ using a global operator and a local oversampling operator, and construct local multi-scale basis functions employing the compactness of the oversampling operator. The resulting method does not involve any global solving of the equation (1.1), and will be detailed in the next section.

## 3 Identify the compact structure of the solution space using Oversampling

In this section, we identify the compact structure of the local solution space through oversampling, and use it to construct basis functions (1.6). In our numerical examples, the domain $\Omega$ is chosen to be $[0,1] \times [0,1]$, and we discretize $\Omega$ using a coarse square mesh of size $H$, which should be chosen according to the desired order of accuracy. With this discetization, we have

$$\Omega = \cup_{i=1}^{N} D_i, \tag{3.1}$$

where $D_i$ have disjoint interiors. Underlying this coarse mesh, we use a triangle fine mesh of size $O(h)$, which is a refinement of the coarse mesh. The fine mesh size $h$ should be chosen such that it can resolve the small scale variation of the multi-scale coefficient in (1.1). In our method we solve the equation (1.1) on the coarse mesh, and the basis functions that we use are constructed and saved using linear basis functions on the fine mesh. The two level discretization is illustrated in Figure 2.

Two Level Mesh

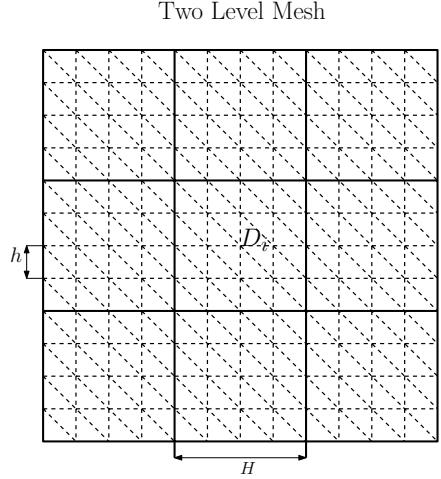

Figure 2

## 3.1 The multi-scale basis

We first introduce a special class of problem-dependent basis functions (1.6), which we call the **multi-scale basis.**

**Definition 3.1** (Multi-Scale basis). For a discretization of $\Omega$ (3.1), we consider basis functions

$$\phi_1(x), \phi_2(x), \ldots \phi_n(x) \in H_0^1(\Omega). \tag{3.2}$$

If they are $a(x)$-harmonic on each coarse element of the coarse discretizaion, $D_j$,

$$-\operatorname{div}(a(x)\nabla \phi_i(x)) = 0, \quad x \in D_j, \tag{3.3}$$

then we call them multi-scale basis functions.

Clearly, the multi-scale basis functions are determined by their restrictions on the boundary of coarse elements $\partial D_i$ since they are $a(x)$-harmonic in each $D_i$. We denote

$$\Gamma = \cup_{i=1}^N \partial D_i. \tag{3.4}$$

We have the following proposition, which implies that if the desired accuracy is $O(H)$, multi-scale basis functions are optimal for fixed local boundary conditions on $\Gamma$ (3.4).

**Proposition 3.1.** *Consider a set of basis functions $\psi_i(x) \in H_0^1(\Omega), i = 1, 2, \ldots m$, and a set of multi-scale basis functions $\phi_i(x) \in H_0^1(\Omega), i = 1, 2 \ldots, n$ on a coarse mesh of size $H$, as showed in Figure 2. Denote the corresponding Galerkin numerical solution (1.8) to (1.1) using $\psi_i(x), i = 1, \ldots m$ as $u_h(x)$, and the Galerkin solution using $\phi_i(x), i = 1, \ldots n$ as $u_h^{MS}(x)$. If*

$$\operatorname{span}\{\phi_1(x)|_\Gamma, \ldots \phi_i(x)|_\Gamma, \ldots \phi_n(x)|_\Gamma\} = \operatorname{span}\{\psi_1(x)|_\Gamma, \ldots \psi_i(x)|_\Gamma \ldots \psi_m(x)|_\Gamma\}. \tag{3.5}$$

*Then we have*

$$\|u(x) - u_h^{MS}(x)\|_a^2 \leq \|u(x) - u_h(x)\|_a^2 + C\|f\|_{L^2(\Omega)}^2 H^2. \tag{3.6}$$

*Namely, if only $O(H)$ accuracy in the energy norm is desired, the multi-scale basis can perform as well as other set of basis functions, given that the local boundary conditions of basis functions are the same.*

To prove the above proposition, we first divide the solution $u(x)$ to (1.1) into two parts. On each coarse mesh element $D_i$, we consider
$$u_i(x) = u(x)|_{D_i}, \tag{3.7}$$
and divide it to an $a(x)$-harmonic part and a local bubble part, as we did in (2.8),
$$u_i(x) = u_i^1(x) + u_i^2(x), \quad x \in D_i, \tag{3.8}$$
where $u_i^1(x)$ is the local $a(x)$-harmonic part, and $u_i^2(x)$ is the local bubble part. Combine these local decompositions from all coarse elements $D_i$ together we get,
$$u(x) = u^1(x) + u^2(x), \quad u^1(x) = \sum_{i=1}^N u_i^1(x), \quad u^2(x) = \sum_{i=1}^N u_i^2(x). \tag{3.9}$$
One can see that the two parts $u^1(x), u^2(x)$ are orthogonal with respect to the $a(\cdot, \cdot)$ inner product (1.9),
$$a(u^1(x), u^2(x)) = 0. \tag{3.10}$$
Besides, the combination of the local bubble parts is small according to (2.12). To be specific,
$$\|u^2(x)\|_a \leq CH \|f\|_{L^2(\Omega)}. \tag{3.11}$$

Next we prove the proposition 3.1.

*Proof.* Denote the numerical solution using $\psi_i(x), i = 1, \ldots m$, $u_h(x)$ as $u_h(x) = \sum_{i=1}^m d_i \psi_i(x)$, then, there exist $c_i, i = 1, \ldots n$ such that $u_h^{ms} = \sum_{i=1}^n c_i \phi_i(x)$, and
$$u_h^{ms}(x)|_\Gamma = u_h(x)|_\Gamma. \tag{3.12}$$
Then we consider
$$\|u_h^{ms}(x) - u(x)\|_a^2 = \|u^2(x) + u^1(x) - u_h^{ms}(x)\|_a^2. \tag{3.13}$$
Since $u_h^{ms}(x) \in H_0^1(\Omega)$ is $a(x)$-harmonic on each coarse element $D_i$, we have
$$a(u^2(x), u_h^{ms}(x)) = \sum_{i=1}^N \int_{D_i} \nabla u_i^2(x)^t a(x) \nabla u_h^{ms}(x) \mathrm{d}x = -\int_{D_i} u_i^2(x) \mathrm{div}(a(x) \nabla u_h^{ms}(x)) \mathrm{d}x = 0. \tag{3.14}$$
Thus $u^2(x)$ is $a$-orthogonal to $u^1(x) - u_h^{ms}(x)$, and according to (3.11) we have
$$\|u(x) - u_h^{ms}(x)\|_a^2 = \|u^2(x)\|_a^2 + \|u^1(x) - u_h^{ms}(x)\|_a^2 \leq \|u^1(x) - u_h^{ms}(x)\|_a^2 + C\|f\|_{L^2(\Omega)}^2 H^2. \tag{3.15}$$
Then we consider $u_e(x) = u(x) - u_h(x)$, and divide it to two parts as we did for $u(x)$ in (3.9). We get
$$u_e(x) = u_e^1(x) + u_e^2(x), \quad a(u_e^1(x), u_e^2(x)) = 0. \tag{3.16}$$
Consequently, we have
$$\|u(x) - u_h(x)\|_a^2 = \|u_e^1(x)\|_a^2 + \|u_e^2(x)\|_a^2 \geq \|u_e^1(x)\|_a^2. \tag{3.17}$$
According to (3.12), we have
$$u_e^1(x) = u^1(x) - u_h^{ms}(x), \tag{3.18}$$
since they are equal on $\Gamma$ and $a(x)$-harmonic on each $D_i$.

Finally based on (3.15), (3.17), and the optimal property (1.10), we have
$$\|u(x) - u_h^{MS}(x)\|_a^2 \leq \|u(x) - u_h^{ms}(x)\|_a^2 \leq \|u(x) - u_h(x)\|_a^2 + C\|f\|_{L^2(\Omega)}^2 H^2, \tag{3.19}$$
and complete the proof. □

As we have showed in (3.11), the bubble part of the solution $u^2(x)$ is small and of $O(H)$ in the energy norm, thus can be neglected if the desired accuracy in the numerical solution is $O(H)$. In our method in this work, we use multi-scale basis functions in (1.6) to approximate the solution space, which are locally $a(x)$-harmonic functions, and are $a(x)$-orthogonal to the bubble part of solution. Due to this $a(x)$-orthogonality and the Galerkin projection formulation in (1.8), multi-scale basis functions only approximate the $a(x)$-harmonic part of the solution and will not bring in additional errors in the bubble part. Thus we can recover the bubble part of solution $u^2(x)$ independently by solve some local bubble problems (2.10). And adding $u^2(x)$ back to $u_h^{MS}(x)$, we can get numerical solution that is free of error in the bubble part. This is one of the advantages of using multi-scale basis functions in (1.6).

To construct local multi-scale basis functions, we divide the $a(x)$-harmonic part of the solution $u^1(x)$ into different parts corresponding to different edges of the coarse mesh, and approximate them separately.

## 3.2 Decomposition of the $a(x)$-harmonic part of the solution

To identify the compact structure of the $a(x)$-harmonic part of the solution, we first introduce a set of primary interpolation multi-scale basis $\psi_i(x), i = 1, \ldots n$ based on the coarse mesh node points $x_1, x_2, \ldots x_n$,

$$\psi_i(x_j) = \delta_{ij}; \quad -\text{div}(a(x)\nabla\psi_i(x)) = 0, \quad x \in D_j. \tag{3.20}$$

We also require that $\psi_i(x)$ is supported on the four coarse elements around $x_i$. For example, we can simply choose the multi-scale basis $\psi_i(x)$ to be linear on the boundaries of coarse elements. We will discuss about the optimal choice of these primary interpolation basis functions in subsection 3.4.

For $f(x) \in L^2(\Omega)$, and the spatial dimension $d = 2, 3$, we have that $u(x)$ is Hölder continuous on $\Omega$ [17], so we can consider the interpolation of $u^1(x)$, namely the $a(x)$-harmonic part of the solution, using the primary basis functions $\psi_j(x)$, and get the residual,

$$u_e^1(x) = u^1(x) - \sum_i u(x_i)\psi_i(x). \tag{3.21}$$

For a coarse mesh element $D_i$, we denote its four nodes points as $x_{i_1}, x_{i_2}, x_{i_3}$ and $x_{i_4}$, then we get the restriction of the residual (3.21) on $D_i$,

$$u_e^1(x)|_{D_i} = u_i^1(x) - u(x_i^1)\psi_{i_1}(x) - u(x_{i_2})\psi_{i_2}(x) - u(x_{i_3})\psi_{i_3}(x) - u(x_{i_4})\psi_{i_4}(x). \tag{3.22}$$

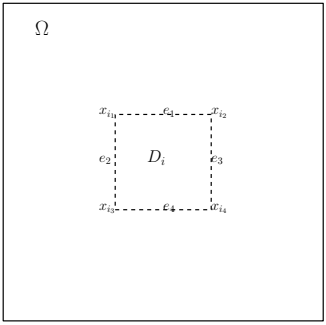
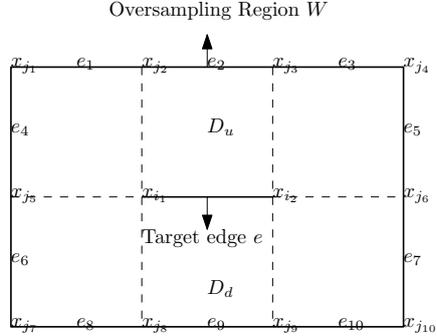

(a) Decomposition of the interpolation residual.  (b) Oversampling region.

Since the residual $u_e^1(x)|_{D_i}$ vanishes on the node points $x_{i_j}, j = 1, 2, 3, 4$, we can divide the trace of $u_e^1$ on $\partial D_i$ to four parts, corresponding to the four edges of $D_i$, $e_1, e_2, e_3, e_4$,

$$u_e^1(x)|_{\partial D_i} = u_e^1(x)|_{e_1} + u_e^1(x)|_{e_2} + u_e^1(x)|_{e_3} + u_e^1(x)|_{e_4}. \tag{3.23}$$

This decomposition is illustrated in Figure 3a. Each part in the above decomposition (3.23) belongs to $H^{1/2}(\partial D_i)$ since they vanish on the node points, and we can extend them to $D_i$ to get four $a(x)$-harmonic components of $u_e^1(x)$. We denote them as $v_{e_1}(x), v_{e_2}(x), v_{e_3}(x), v_{e_4}(x)$, and have

$$u_e^1(x)|_{\partial D_i} = v_{e_1}(x) + v_{e_2}(x) + v_{e_3}(x) + v_{e_4}(x). \tag{3.24}$$

Combining these local decompositions together, we have

$$u_e^1(x) = \sum_e v_e(x), \tag{3.25}$$

where $v_e(x)$ is the $a(x)$-harmonic extension of the interpolation error on the edge $e$ to its two neighbor elements. In the above decomposition (3.25), we are actually dividing the error $u_i^e(x)$ in the $a(x)$-harmonic part of solution into different parts corresponding to errors on different edges $e$. This is possible since $u_i^e(x)$ vanishes on the node points, thanks to the interpolation operation using $\psi_i(x)$ (3.21).

We seek to construct boundary basis functions on each edge $e$ that approximate $v_e(x)$, and combine them together to get the whole trial space. We introduce the following operator for the edge $e$ with endpoints $x_{i_1}$ and $x_{i_2}$, which maps $f(x) \in L^2(\Omega)$ to the interpolation residual of the solution on $e$,

$$T_e: \quad f(x) \in L^2(\Omega) \to v_e(x) = u(x) - u(x_{i_1})\phi_{i_1}(x) - u(x_{i_2})\phi_{i_2}(x) \in H^{1/2}(e). \tag{3.26}$$

The left singular vectors of $T_e$ form the optimal local boundary basis functions. However, $T_e$ is a global operator and its construction involves solving the equation (1.1) globally. In the next section, we decompose $T_e$ as a global solution operator and a local oversampling operator, and construct boundary basis functions that approximate the range of $T_e$ through the oversampling operator.

## 3.3 The oversampling operator

To identify the compact structure of the solution space restricted on the edge $e$, we put it in an oversampling region that we denote by $W$. In our numerical examples, we use the square mesh of size $H$ for the coarse discretization, and the oversampling region $W$ is chosen as the union of the six elements around the edge $e$. It is illustrated in Figure 3b. We remark that our method is also applicable to other types of discretizations like triangular mesh. We denote the solution on $W$ as $u_W(x)$.

We remark that the idea of identifying the local structure of the solution space by putting it in a larger region, namely oversampling, was first proposed in [18] to reduce the resonance error due to artificial local boundary conditions, and this strategy was later employed in [2, 5, 10].

We denote $T_W$ as the operator that maps $f(x) \in L^2(\Omega)$ to the oversampling solution $u_W(x) = u(x)|_W$, and $T_{W \to e}$ as the operator that maps $u_W(x)$ to the solution restricted on the edge $e$:

$$T_W: \ f(x) \to u_W(x) = u(x)|_W, \quad T_{W \to e}: \ u_W(x) \to u_e(x) = u_W(x)|_e. \tag{3.27}$$

We also introduce the interpolation residual operator using (3.20), $P_e$,

$$P_e: \ u(x)|_e \to u(x)|_e - u(x_{i_1})\psi_{i_1}(x) - u(x_{i_2})\psi_{i_2}(x). \tag{3.28}$$

With the above definitions, the operator $T_e$ (3.26) can be decomposed as

$$T_e = P_e T_{W \to e} T_W. \tag{3.29}$$

We call the operator $P_e T_{W \to e}$ in the above decomposition (3.29) the oversampling operator, which maps the solution on $W$, $u_W(x)$ to the interpolation residual, which we denote by $v_e(x)$,

$$P_{OS} = P_e T_{W \to e}: \quad u_W(x) \to v_e(x) = u_W(x) - u_W(x_{i_1})\psi_{i_1}(x) - u_W(x_{i_2})\psi_{i_2}(x), \tag{3.30}$$

where $x_{i_1}$ and $x_{i_2}$ are the two endpoints of $e$, and $\psi_i(x)$ is the primary interpolation basis (3.20).

The solution to (1.1) on $W$, $u_W(x)$ can be divided into two parts, the $a(x)$-harmonic part $u_W^1(x)$ and the bubble part $u_W^2(x)$. We employ the compactness of the oversampling operator (3.30) to construct basis functions in $H^{1/2}(e)$ that vanish at $x_{i_1}$ and $x_{i_2}$, and approximate the range of (3.26). To be specific, we use the first several left singular vectors of $P_{OS}$ as the basis functions associated with $e$. We first introduce appropriate inner products for the domain and range space of $P_{OS}$.

On the edge $e$, the image of $T_e$, $v_e(x) \in H^{1/2}(e)$ and vanishes on the two endpoints. We consider its $a(x)$-harmonic extension to the upper and lower coarse elements respectively, as showed in Figure 3b, and denote them as $v_e^u(x)$ and $v_e^d(x)$. Then we define

$$\|v_e(x)\|^2_{H^{1/2}(e)} = \frac{1}{2}\int_{D_u} \nabla v_e^u(x)^t a(x) \nabla v_e^u(x) \mathrm{d}x + \frac{1}{2}\int_{D_d} \nabla v_e^d(x)^t a(x) \nabla v_e^d(x) \mathrm{d}x. \tag{3.31}$$

In the domain of the operator $P_e T_{W \to e}$, namely, $u_W(x)$, we define its inner product as

$$\|u_W(x)\|_a^2 = \int_W \nabla u_W^1(x)^t a(x) \nabla u_W^1(x) + (u_W^1(x))^2 \mathrm{d}x + \int_W [\mathrm{div}(a(x) \nabla u_W(x))]^2. \tag{3.32}$$

With the above inner products, we compute the singular value decomposition of the oversampling operator $P_{OS}$. To discretize the domain of $P_{OS}$, we consider its two parts, the $a(x)$-harmonic part, and the bubble part. The $a(x)$-harmonic part only depends on the trace of $u_W(x)$ on $\partial W$, we discretize $H^{1/2}(\partial W)$ using all the fine mesh piecewise linear functions. If $\partial W$ intersects with $\partial \Omega$, then we choose $H^{1/2}(\partial W)$ to be fine mesh basis functions that vanish on $\partial \Omega$. The bubble part of the solution $u_W^2(x)$ only depends on $f_W(x) = f(x)|_W$, and we discretize $f(x)$ using piecewise constant functions on the coarse mesh. The following lemma justifies this discretization of $f(x) \in L^2(\Omega)$.

**Lemma 3.1.** *Denoe the space of piecewise constant functions on the coarse mesh as $V_c$. The coarse mesh is showed in Figure 2 with mesh size $H$. Then for $f(x) \in L^2(\Omega)$, there exist $f_c(x) \in V_c$, such that*

$$\|f(x) - f_c(x)\|_{H^{-1}(\Omega)} \leq CH \|f(x)\|_{L^2(\Omega)}. \tag{3.33}$$

According to (1.3), the above Lemma implies that if we want to obtain $O(H)$ accuracy in the energy norm of the numerical solution, we can only consider piecewise constant forcing functions on the coarse mesh. If higher accuracy is desired in the numerical solution for $f(x) \in L^2(\Omega)$, we can simply refine the discretization of $L^2(\Omega)$ in the oversampling operator (3.30). We remark that if $f(x)$ has higher regularity, we can get better approximation result in (3.33), for example

$$\inf_{f_c \in V_c} \|f(x) - f_c(x)\|_{H^{-1}(\Omega)} \leq CH^2 \|f(x)\|_{H^1(\Omega)}. \tag{3.34}$$

With the above discretization of $u_W(x)$, we truncate the singular values of $P_{OS}$ to $O(H)$ and select the corresponding left singular vectors as the boundary basis functions $e$. We denote them as

$$v_e^1(x), \ldots v_e^{k_e}(x) \in H^{1/2}(e), \tag{3.35}$$

which vanish on the two endpoints of $e$. According to the oversampling operator in (3.29), and our truncation criteria, we have the following approximation property

$$\inf_{c_i^e} \|u(x)|_e - u(x_{i_1})\phi_{i_1}(x) - u(x_{i_2})\phi_{i_2}(x) - \sum_{i=1}^{e_k} c_i^e v_e^i(x)\|_{H^{1/2}(e)} \leq CH(\|u_W^1(x)\|_a + \|f\|_{L^2(W)}). \tag{3.36}$$

Then we extend these boundary basis functions to the two neighbourhood coarse elements $D_u$ and $D_d$ as $a(x)$-harmonic functions, by solving

$$\begin{cases} -\text{div}(a(x)\nabla \phi_e^k(x)) = 0, x \in D_u, D_d, \\ \phi_e^k(x)|_e = v_e^k(x). \end{cases} \tag{3.37}$$

Finally, we combine the multi-scale basis functions associated with each edge of the coarse mesh, $e$,

$$\phi_e^i(x), \quad i = 1, 2 \ldots k_e. \tag{3.38}$$

with the primary interpolation basis functions $\psi_i(x)$ (3.20), and get the trial space,

$$V_h = \text{span}\{\phi_i(x), \quad i = 1, \ldots n\}. \tag{3.39}$$

We have the following error estimate using the trial space (3.39).

**Proposition 3.2.** *Using the trial space consisting of the primary interpolation basis functions (3.20) and the multi-scale basis functions constructed from each edge of the coarse mesh (3.37), we obtain numerical solution to (1.1) using the Galerkin projection (1.8). Then we have the following convergence property,*

$$\|u(x) - u_h^{MS}(x)\|_a \leq CH \|f(x)\|_{L^2}. \tag{3.40}$$

**Remark 3.1.** *Using a simple Aubin-Nitsche duality argument and the convergence result in the energy norm (3.40), we can get the following convergence result in the $L^2$ norm.*

$$\|u(x) - u_h^{MS}(x)\|_{L^2(\Omega)} \leq CH^2 \|f(x)\|_{L^2}. \tag{3.41}$$

The proof of the convergence result (3.40) follows directly from the decomposition of the solution operator (3.29) and the truncation in the singular value decomposition of the oversampling operator.

*Proof.* We choose $c_e^j$ as the ones in (3.36), and denote

$$u_h^{ms}(x) = \sum_{i=1}^n u(x_i)\psi_i(x) + \sum_e \sum_{j=1}^{k_e} c_e^j \phi_e^j(x) \in V_h. \tag{3.42}$$

Then we consider $\|u_h^{ms}(x) - u(x)\|_a$, since the basis functions in (3.39) are multi-scale basis, namely, they are $a(x)$-harmonic on each $D_i$, we have

$$\|u_h^{ms}(x) - u(x)\|_a^2 = \|u_h^{ms}(x) - u_1(x) - u_2(x)\|_a^2 \leq \|u_h^{ms}(x) - u_1(x)\|_a^2 + CH^2 \|f\|_{L^2(\Omega)}^2, \tag{3.43}$$

where $u_1(x)$ and $u_2(x)$ are the $a(x)$-harmonic part and bubble part of the solution $u(x)$.

Then we divide $\|u_h^{ms}(x) - u_1(x)\|_a^2$ to different parts on $D_i$. For each part, according to the approximation property (3.36), and the definition (3.31), we have

$$\int_{D_i} \nabla(u_h^{ms}(x) - u_1(x))^t a(x) \nabla(u_h^{ms}(x) - u_1(x)) \text{d}x \leq C \sum_W H \|u_W^1(x)\|_a^2. \tag{3.44}$$

The sum over $W$ corresponds to the oversampling regions for edges of $D_i$. There are four of them for each $D_i$. Summing up (3.44) for all the coarse elements of $\Omega$, we have

$$\|u_h^{ms}(x) - u_1(x)\|_a^2 \leq CH^2 \|u(x)\|_a. \tag{3.45}$$

Putting (3.45) in (3.43), we have

$$\|u(x) - u_h^{ms}(x)\|_a \leq CH\|f(x)\|_{L^2(\Omega)}. \tag{3.46}$$

Then using the optimal approximation property (1.10), we finish the proof. $\square$

To make the number of multi-scale basis functions in (1.6) small, we want the singular values of $P_{OS} = P_e T_{W \to e}$ decay fast. $P_{OS}$ can be divided into two parts: the first part acts on the $a(x)$-harmonic part of $u_W(x)$, and we denote it as $P_{OS}^1$; the second part acts on the bubble that only depends on $f_W(x)$, and we denote it as $P_{OS}^2$. For the first part, similar analysis has been done in [5] in a slightly different setting, and the method there also applies to our problem. We have the following result.

**Proposition 3.3.** *Denote the singular values of $P_{OS}^1$ that acts on $a(x)$-harmonic functions on $W$ as $\sigma_k$, then we have the following upper bound on the decay of $\sigma_k$: for any $\epsilon > 0$, there exist $C$ such that*

$$\sigma_k \leq C \exp\{-k^{1/(d+1)-\epsilon}\}, \tag{3.47}$$

*where $d$ is the dimension of the domain $\Omega$.*

The second part $P_{OS}^2$ is small according to (2.12), and the decay rate of its singular values can be obtained from (2.5) using a simple scaling argument.

We will see in our numerical results section that the singular values of $P_{OS}$ decay very fast, and a very small number of boundary basis functions can achieve high local approximation accuracy.

As we have shown previously, the basis functions we obtain are multi-scale basis functions, thus are $a(x)$-orthogonal to the bubble part of the solution space. This gives us the flexibility to add the bubble parts back to the numerical solutions at local regions where higher accuracy is desired by simply solving some local bubble problems (2.10). In our truncation of the singular values of the local compact operator $P_e T_{W \to e}$, we choose the threshold to be $O(H)$, since $O(H)$ accuracy is required in (3.40). If we need higher accuracy than $H$, for example, $O(\epsilon)$ where

$$h \ll \epsilon \ll H. \tag{3.48}$$

Then we can truncate the singular values of $P_e T_{W \to e}$ by $\epsilon$, by doing which, the resulting multi-scale basis functions are able to approximate the $a(x)$-harmonic part of the solution space up to $O(\epsilon)$ accuracy. Then by adding back the bubble part of the solution $u^2(x)$ to the numerical solution $u_h^{MS}(x)$, we can get $O(\epsilon)$ accuracy in our final numerical solutions. Namely, our upscaling strategy allows us to get arbitrarily high accuracy that is permitted by the fine mesh discretization.

### 3.4 Optimal primary interpolation basis functions

In constructing the multi-scale basis functions in the previous section, we need to choose a set of primary interpolation basis functions (3.20) first, which allows us to reduce the problem to approximating the solutions space restricted on each edge $e$, (3.29). The choice of these interpolation basis will affect the oversampling operator $P_{OS} = P_e T_{W \to e}$. In this subsection, we identify the optimal interpolation basis functions by solving local under-determined least square problems.

The oversampling operator for edge $e$, $P_{OS}$, depends on the interpolation basis functions $\psi_{i_1}(x)$ and $\psi_{i_2}(x)$ associated with the two endpoints of $e$, and we seek optimal primary interpolation basis functions $\psi_{i_1}(x)$ and $\psi_{i_2}(x)$ that make the singular values of $P_{OS}$ have the fastest decay.

We consider the following optimization problem

$$\min_{\phi_{i_1}(x)} \|\psi_{i_1}(x)\|_a, \quad \min_{\phi_{i_2}(x)} \|\psi_{i_2}(x)\|_a, \quad \text{subject to} \tag{3.49a}$$

$$-\text{div}(a(x)\nabla \psi_{i_j}(x)) \in L^2(W), \; x \in W, \quad \psi_{i_j}(x_{i_k}) = \delta_{jk}. \tag{3.49b}$$

where the norm $\|\cdot\|_a$ is defined in (3.32).

Under certain conditions, the optimization problem (3.49) has unique solutions, and the unique solutions are optimal interpolation basis functions. We have the following theorem.

**Theorem 3.1.** *Let $\phi_{i_1}^*(x), \phi_{i_2}^*(x)$ be the solution to (3.49), and $\phi_{i_1}(x)$ and $\phi_{i_2}(x)$ be two other multi-scale interpolation basis functions. Denote the corresponding oversampling operators using these interpolation basis as (3.30) as $P_{OS}^*$ and $P_{OS}$, and let $\sigma_i^*, i = 1, 2, \ldots$ and $\sigma_i, i = 1, 2, \ldots$ be their singular values. Assume that both $\sigma_i^*$ and $\sigma_i$ are arranged in descending order. Then we have*

$$\sigma_i^* \leq \sigma_i, \tag{3.50}$$

*which implies the interpolation basis in (3.49) make the singular value of the oversampling operator (3.30) have the fastest decay.*

The proof of Theorem 3.1 is given in the appendix.

Note that in the minimization problem (3.49), the two interpolation basis functions can be constructed independently by solving under-determined least square problems. By going over the oversample regions for each edge of the coarse mesh, we can construct the optimal interpolation basis functions (3.20) on all the boundaries of local regions. Then we can extend them locally to $a(x)$-harmonic functions by solving some local boundary value problems as in (3.37) to get the interpolation basis (3.20).

We will see in our numerical results section that, in some cases, the optimal interpolation basis functions (3.49) themselves are enough to approximate the solution space of (1.1). Namely, there is no need to construct the multi-scale basis functions associated with each edge (3.37).

We also remark that the optimization problem (3.49) can be interpreted as a Bayesian Inference problem as in [29, 26, 25]. Interested readers may consult the references for more information.

### 3.5 Implementation of the whole method

The proposed multi-scale finite element method consists of two stages, the offline stage and the online stage. In the offline stage, we identify the compact structure of the solution space. In the online stage, for a given forcing function, we compute the numerical solution using the offline basis functions.

The offline stage involves the following procedures,

1. Build the $a(x)$-harmonic extension operator on each oversampling region $W$.

   On each oversampling region $W$, we build the local $a(x)$-harmonic extension operator, which maps the boundary condition which belongs to $H^{1/2}(\partial W)$ to $a(x)$-harmonic functions on $W$. This step requires solving a series of boundary value problems on $W$.

2. Construct the $H^{1/2}(\partial W)$ norm.

   The $a(x)$-harmonic extension operator and the inner product (3.32) together induces a norm on $H^{1/2}(W)$. We need this weighted norm in the SVD of the oversampling operator.

3. Construct the optimal multi-scale interpolation basis functions, $\psi_{i_1}(x)$ and $\psi_{i_2}(x)$.

   Discretize the domain of the oversampling operator $P_{OS}$ using the strategy described in section 3.3, and build the restriction operator that maps $u_W(x)$ to $u_W(x_{i_1})$ and $u_W(x_{i_2})$. Then we solve the optimization problem (3.49), which is essentially an under-determined least square problem.

   Solving the optimization problem (3.49) on each oversampling region only gives us the boundary conditions of the primary interpolation basis functions, we need to solve local boundary value problems to get their local $a(x)$-harmonic extensions as the interpolation basis functions.

4. Build the oversampling operator.

   Using the optimal interpolation basis functions we constructed in the previous step, we build the oversampling operator $P_{OS}$, (3.3).

5. Construct the $H^{1/2}(e)$ norm.

   Construct the $H^{1/2}(e)$ norm (3.31), which requires solving a series of local boundary value problems on the two neighbour coarse elements of $e$.

6. Compute the SVD of the oversampling operators.

   Using the inner product (3.32) and (3.31) in the oversampling operator to compute its singular value decomposition. We truncate the singular values to $\epsilon$, and save the corresponding left singular vectors, which are basis functions of $H^{1/2}(e)$, $v_e^1(x), v_e^2(x), \ldots v_e^{k_e}(x)$.

7. Construct multi-scale basis functions associated with each edge $e$.

   We extend the basis functions of $H^{1/2}(e)$, namely $v_e^i(x)$, to the two neighbourhood coarse elements of $e$ by solving local boundary value problems (3.37).

8. Compute the stiffness matrix.

   Combining the multi-scale interpolation basis functions with the multi-scale basis functions associated with each edge of the coarse mesh, we get the trial space

   $$V_h = \text{span}\{\phi_1(x), \phi_2(x), \ldots \phi_n(x)\}. \tag{3.51}$$

   Then we save these multi-scale basis functions, and compute the stiffness matrix,

   $$M(i,j) = a(\phi_i(x), \phi_j(x)). \tag{3.52}$$

The online stage involves the following procedures,

1. Compute the load vector.

   For a given forcing function $f(x) \in L^2(\Omega)$, we compute the corresponding load vector

   $$b(i) = \int_\Omega \phi_i(x) f(x) \mathrm{d}x, \quad i = 1, \ldots n. \tag{3.53}$$

2. Compute the online numerical solution.

   Using the load vector $b$ (3.53) and stiffness matrix (3.52), we solve the linear system

   $$Mc = b, \tag{3.54}$$

   and get the online numerical solution on the fine mesh

   $$u_h^{MS}(x) = \sum_{i=1}^n c_i \phi_i(x). \tag{3.55}$$

   Recall that the multi-scale basis functions associated with the edges vanish on the coarse grid node points. So if we only want the large scale solution, we can simply select the coefficients in $c$ that correspond to the primary multi-scale interpolation basis (3.20).

3. Recover the bubble part of the solution.

   Solve the local boundary value problem (2.10) on each coarse mesh element $D_i$, and get $u_i^2(x)$. Combine these local bubble parts together and add them back to Galerkin solution (3.55), we get

   $$u_h(x) = u_h^{MS}(x) + \sum_{i=1}^N u_i^2(x). \tag{3.56}$$

   If $O(H)$ accuracy is required in the numerical solution, this step is unnecessary since the bubble part does not impact the large scale properties of the solution.

Note that in the offline stage, we need to solve a series of boundary value problems for each edge of the coarse mesh to construct the oversampling operator (3.30), and then compute its singular value decomposition, which is computationally expensive. However, the constructions of multi-scale basis functions on each edge are independent from each other, thus the offline stage can be implemented on a parallel machine to accelerate the computation. In the online stage, the main computational cost comes from solving the linear system (3.54). Our numerical results in the next section suggest that a small number of basis functions are enough to obtain the coarse mesh accuracy, $O(H)$, thus the linear system (3.54) is small and sparse. This implies that the online computation in our method is efficient, and our method can bring in significant computational savings in the multi-query setting, where the equation (1.1) needs to be solved for multiple times using different forcing functions.

# 4  Numerical Results

In this section, we present numerical examples that have multiple-scale features and high-contrast channels to demonstrate the capacity of our method in identifying and exploiting the compact structure of the local solution space to achieve computational savings (in the online stage). We discretize the domain of the problems $\Omega = [0,1] \times [0,1]$ using a two level mesh as showed in Figure 2. The coarse mesh is of size $H = 1/32$, and the fine mesh is of size $h = 1/1024$.

## 4.1 An example with multiple spatial scales

The first example we consider is one that has multiple spatial scales. The coefficient is given by (2.14), and it is visualized in Figure 1a. For each edge of the coarse mesh, we compute the singular value decomposition of the oversampling operator, and truncate the singular values at $\epsilon = H$, which guarantees $O(H)$ accuracy in the online numerical solution. After the offline stage, multi-scale basis functions (1.6) are constructed, and the average number of basis functions associated with each edge is

$$\bar{k}_e = \frac{\sum_e k_e}{\#(e)} \approx 1.00. \tag{4.1}$$

$\bar{k}_e$ is very small. Actually only 1 or 2 multi-scale basis functions are constructed for each edge of the coarse mesh, aside from the primary interpolation basis functions. And this reflects the efficiency of our method in the online stage since the stiffness matrix is small and sparse.

To measure the error in our online numerical solution, we need to choose a reference solution. Since the multi-scale basis functions are constructed and saved on the fine mesh of size $h$, we will use the piecewise linear finite element solution on the fine mesh as the reference.

In the online stage, we choose the forcing function $f(x, y)$ to be

$$f(x, y) = 1, \quad (x, y) \in \Omega. \tag{4.2}$$

Recall that the basis functions that we use are $a(x)$-harmonic in each $D_i$, and we can add back the bubble part of the solution by simply solving some local cell problems. Our online numerical solution, and the numerical solution after correction using the bubble part are plotted in Figure 4. We can see that the numerical errors in the online numerical solution is very small.

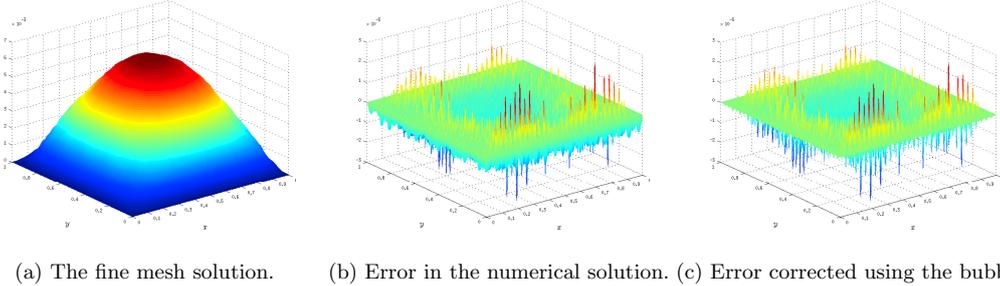

(a) The fine mesh solution.    (b) Error in the numerical solution. (c) Error corrected using the bubble.

Figure 4: Online numerical solutions

We measure the error of the numerical solution in the energy norm and $L^2$ norm. We denote the numerical solution as $u_H^{MS}(x)$, and the corrected solution using the bubble part as $u_H(x)$. We compute

$$E_a^{MS} = \frac{\|u(x) - u_H^{MS}(x)\|_a}{\|u(x)\|_a}, \qquad E_a = \frac{\|u(x) - u_H(x)\|_a}{\|u(x)\|_a}, \tag{4.3a}$$

$$E_{L^2}^{MS} = \frac{\|u(x) - u_H^{MS}(x)\|_{L^2(\Omega)}}{\|u(x)\|_{L^2(\Omega)}}, \qquad E_{L^2} = \frac{\|u(x) - u_H(x)\|_{L^2(\Omega)}}{\|u(x)\|_{L^2(\Omega)}}. \tag{4.3b}$$

The results are listed in Table 1. We can see that by adding the bubble part of solution back to the numerical solution, the error in $L^2$ norm and energy norm are both reduced by about one half. This implies the numerical error in the $a(x)$-harmonic part of the solution is about the same as that in the bubble part. The latter is of order $O(H)$ in the energy norm and of order $O(H^2)$ in the $L^2$ norm, and this result confirms our error estimates (3.40) and (3.41).

Then we consider only using the primary interpolation basis functions $\psi_i(x), i = 1, \ldots N$ in the trial space (3.51), namely, we do not use the multi-scale basis functions associated with each edge of the coarse mesh. The error in the corresponding numerical solution is also listed in Table 1. We can that the relative error in $L^2$ is also small, which means the numerical solution can capture the large-scale property of the solution. However, the errors in the energy norm and $L^2$ norm are both significantly larger than that in (4.3), which implies the necessity of enriching the trial space using the multi-scale basis functions associated with the edges of the coarse mesh if higher accuracy is required.

| | Energy Norm | $L^2$ Norm |
|:---:|:---:|:---:|
| Numerical solution | $4.16 \times 10^{-2}$ | $1.73 \times 10^{-3}$ |
| Corrected numerical solution | $2.67 \times 10^{-2}$ | $8.75 \times 10^{-4}$ |
| Interpolation basis only | $7.57 \times 10^{-2}$ | $5.95 \times 10^{-3}$ |

Table 1: Numerical errors

To further demonstrate the convergence rate in (3.40) and (3.41), we consider a sequence of coarse mesh size with $H = 2^{-k}, k = 3, 4, 5, 6, 7$. For each $H$, we compute the error in the online numerical solution, the decay of the numerical error with $H$ is plotted in Figure 5. From the plot, we can clearly see that the decay rates of the online numerical error agree with the estimates (3.40) and (3.41).

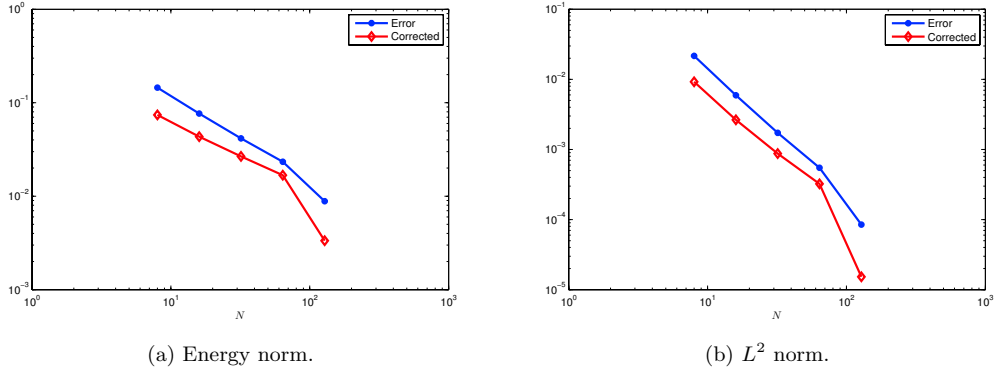

(a) Energy norm.  (b) $L^2$ norm.

Figure 5: Convergence of the online numerical solution with $N = 1/H$.

### 4.2 An example without scale-separation

In this subsection, we consider an example of the coefficient $a(x)$ without scale-separation.

$$a(x, y) = |\tilde{a}| + 0.5, \qquad (4.4)$$

where $\tilde{a}$ is normally distributed on the node point of a mesh of size $\frac{1}{128}$. And $a(x, y)$ on the fine mesh is obtained using piecewise linear interpolation. This coefficient $a(x, y)$ is rough and has no clear scale-separation. The configuration of the coefficient (4.4) is illustrated in Figure 6.

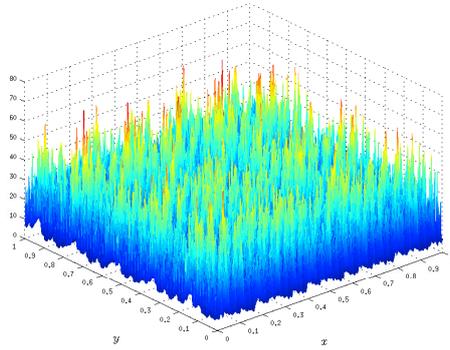

Figure 6: The rough coefficient $a(x)$ without scale-separation.

We discretize the spatial domain $\Omega$ using a two level mesh as showed in Figure 2, and then we solve the optimization problem (3.49) and build the oversampling operator (3.30) on the fine mesh. We truncate the singular value decomposition of the oversampling operator to $\epsilon = H$. After the offline stage, the average number of multi-scale basis functions associated with each edge of the coarse mesh, (4.1), is $\bar{k}_e \approx 1.00$. The smallness of $\bar{k}_e$ reflects the compactness of the local solution space, and the efficiency of our method in the offline stage since the stiffness matrix is small and sparse.

|                              | Energy Norm           | $L^2$ Norm            |
|------------------------------|-----------------------|-----------------------|
| Numerical solution           | $4.44 \times 10^{-2}$ | $1.98 \times 10^{-3}$ |
| Corrected numerical solution | $3.17 \times 10^{-2}$ | $1.18 \times 10^{-3}$ |

Table 2: Numerical errors

In the online stage, we choose the forcing function $f(x)$ to be same as (4.2), and measure the error of the online numerical solution using the fine mesh solution as reference. The numerical solutions are plotted in Figure 7. We can see that the errors in the online numerical solutions are small. We measure the error (4.3) in the energy norm and the $L^2$ norm. The results are summarized in Table 2. Again we see that our method achieves very high accuracy in the online stage, which reflects that the good performance of our method does not depend on the scale-separation of the coefficient.

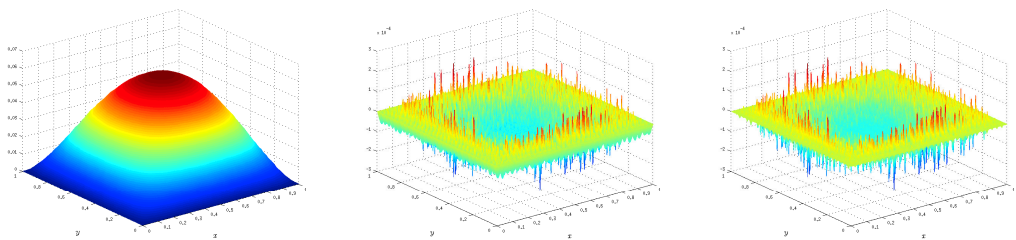

(a) The fine mesh solution.    (b) Error in the numerical solution. (c) Error corrected using the bubble.

Figure 7: Online numerical solutions

## 4.3 An example with high-contrast channels

In this subsection, we consider an example with high-contrast channels. The high contrast in the coefficient violates our uniform ellipticity assumption (1.2), and brings in additional difficulty. The coefficient that we consider here is the one with multiple scales (2.14) added with some high conductivity patches and channels. $\log_{10} a(x)$ is plotted in Figure 8, which has very strong heterogeneity.

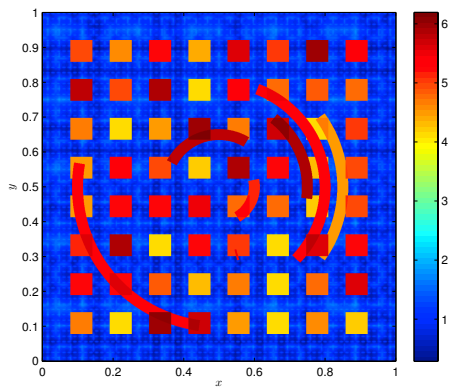

Figure 8: $\log_{10} a(x)$. The coefficient with high-contrast channels.

We discretize the problem in the spatial domain as the previous two examples, and build the oversampling operator for each edge $e$ of the coarse mesh. We truncate the singular value decomposition of the oversampling operators to $\epsilon = H$. The average number of multi-scale basis functions associated with each edge is $\bar{k}_e = 0.89$ (4.1). Namely, on average, we use less than one multi-scale basis function for each edge of the coarse mesh, which implies the compactness of the solution space on local regions of the domain for this problem with high contrast coefficient.

|                              | Energy Norm          | $L^2$ Norm           |
|------------------------------|----------------------|----------------------|
| Numerical solution           | $3.67 \times 10^{-2}$ | $1.64 \times 10^{-3}$ |
| Corrected numerical solution | $2.02 \times 10^{-2}$ | $6.13 \times 10^{-4}$ |

Table 3: Numerical errors

In the online stage, we choose the forcing function to be (4.2), and use fine mesh solution as reference. The numerical errors are plotted in Figure 9. We can see that our numerical solutions have high accuracy and can capture the large-scale properties of the solution. The numerical errors of the solutions (4.3)

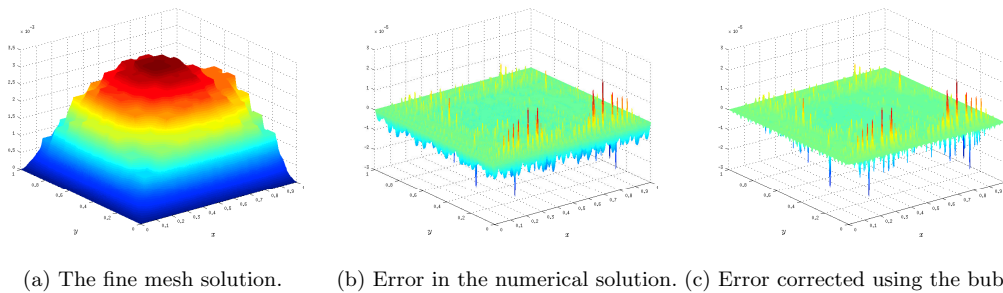

(a) The fine mesh solution.   (b) Error in the numerical solution.   (c) Error corrected using the bubble.

Figure 9: Online numerical solutions.

are listed in Table 3. Again, we see the high accuracy of our online numerical solution.

## 5 Concluding Remarks

In this paper, a novel multi-scale finite-element method is proposed, which is based on the compactness of the solution space on local regions of the spatial domain, and does not depend on any scale-separation or periodicity assumption of the coefficient. By introducing a set of primary interpolation basis functions, we reduce the problem of approximating the local solution space to approximating the trace of the solution on each edge of the coarse mesh. We construct basis functions for each edge of the coarse mesh separately employing an oversampling operator which is local and compact. The optimal primary interpolation basis is also identified by solving a series of under-determined least square problems.

The resulting method involves two stages: in the offline stage we identify the local compact structure of the solution space; in the offline stage we solve the problem efficiently exploiting this compact structure. The offline computation is expensive but it only requires solving local problems and is parallel in nature. We can rigorously control the error in the online numerical solutions by thresholding in the offline stage. The resulting basis functions are called multi-scale basis because they are $a(x)$-harmonic on each coarse element, which have optimal approximation property. Multi-scale basis functions are orthogonal to the bubble part of the solution space, and this gives us the flexibility to put back the bubble part of the solution in regions where high accuracy is desired.

Numerical results suggest that our method is very robust and can achieve high accuracy for the challenging problems without scale-seperation, or have high-contrast inclusions.

## Acknowledgements

This research was in part supported by Air force MURI Grant FA9550-09-1-0613, DOE Grant DE-FG02-06ER257, and NSF Grant No. DMS-1318377, DMS-1159138.

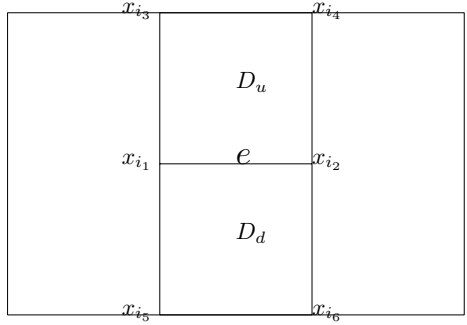

Figure 10: The oversampling region.

## APPENDIX

# Proof of Theorem 3.1

Recall that the restriction operator $P_{W \to e}$ (3.27), maps $u(x)|_W$ to $H^{1/2}(e) \cap C^\alpha(e)$,

$$P_{W \to e}: \ u_W(x) \to u_e(x) \in H^{1/2}(e) \cap C^\alpha(e). \tag{.1}$$

And we have introduced the following inner product (3.32) on the domain of $T_{W \to e}$, which denote by $V_W = D(T_{W \to e})$,

$$\|u_W(x)\|_{V_W}^2 = a(u_W^1(x), u_W^1(x)) + \|u_W^1(x)\|_{L^2(W)}^2 + \|\text{div}(a(x)\nabla u_W(x))\|_{L^2(W)}. \tag{.2}$$

Note that $u_e(x)$ does not necessarily vanish on the two endpoints of $e$, thus the inner product we defined in (3.31) does not apply to $v_e(x)$. We first extend $u_e(x)$ to the boundary of its two neighborhood coarse elements $D_u$, $D_d$, Figure 10. To be specific, we let $u_e(x)$ vanish on $x_{i_3}, x_{i_4}, x_{i_5}, x_{i_6}$, and be piecewise linear on $\partial(\bar{D}_u \cup \bar{D}_d)$. Then we solve local boundary value problems on $D_u$ and $D_d$ separately to extend $u_e(x)$ to $D_u \cup D_d$ as $a(x)$-harmonic functions. With this $a(x)$-harmonic extension to $D_u$ and $D_d$, we introduce the following inner product, which agrees with (3.31) for $u_e(x)$ that vanishes on $x_{i_1}$ and $x_{i_2}$.

$$\|u_e(x)\|_{H^{1/2}(e)}^2 = \int_{D_u \cup D_d} \nabla u_e(x)^t a(x) \nabla u_e(x) \mathrm{d}x. \tag{.3}$$

And (.3) can be viewed as an extension of (3.31).

Next we consider $P_1$ and $P_2$, which are the bounded linear functional that maps $u_W(x) \in V_W$ to its values on $x_{i_1}$ and $x_{i_2}$ separately. We can easily see that $P_1$ and $P_2$ are linearly independent, thus there exist $\phi_{i_1}(x)$ and $\phi_{i_2}(x)$, such that

$$P_j(\phi_{x_{i_k}}(x)) = \delta_{jk}, \quad j, k = 1, 2. \tag{.4}$$

We denote the intersection of the kernel of $P_1(x)$ and $P_2$ as $V_W^0$, which is closed subspace of $V_W$. And we denote the projection of $\psi_{i_1}(x)$ and $\psi_{i_2}(x)$ to the orthogonal complement of $V_W^0$, $(V_W^0)^\perp$, as $\psi_{i_1}^*(x)$ and $\psi_{i_2}^*(x)$, namely,

$$\psi_{i_1}(x) - \psi_{i_1}^*(x), \ \psi_{i_2}(x) - \psi_{i_2}(x) - \psi_{i_2}^*(x) \perp V_W^0, \tag{.5}$$

where the orthogonality is in the sense of (.2). Then we have

$$\psi_{i_j}^*(x_{i_k}) = \delta_{jk}, \quad \|\psi_{i_j}(x)\|_{V_W} \geq \|\psi_{i_j}^*(x)\|_{V_W}, \quad j, k = 1, 2. \tag{.6}$$

And we finish the proof of the existence of unique solutions $\psi_{i_1}^*(x)$ and $\psi_{i_2}(x)$ to the optimization problem (3.49). Next we proof the property (3.50).

We choose $\psi_{i_1}^*(x)|_e$ and $\psi_{i_2}^*(x)|_e$ as the interpolation basis in (3.3), and get the oversampling operator $P_{OS}^*$. For any other two interpolation basis functions, $\psi_{i_1}(x)$ and $\psi_{i_2}(x)$, we denote the corresponding oversampling operator as $P_{OS}$, then we compare their singular values, $\sigma_k^*$, and $\sigma_k$. We use the following characterization of singular values,

$$\sigma_k = \sup_{V_k} \inf_{\|u(x)\|_{V_W}=1} \|P_{OS}u(x)\|_{H^{1/2}(e)}. \tag{.7}$$

where $V_k$ is a $k$-dimensional subspace of $V_W$.

For any $u(x) \in V_W$, we consider its decomposition in $V_W^0$ and $(V_W^0)^\perp$, and denote them as $u_1(x)$, $u_2(x)$,

$$u(x) = u_1(x) + u_2(x). \tag{.8}$$

Then according to our choice the interpolation basis function, we have

$$P_{OS}^*(u_1(x)) = 0, \quad P_{OS}^*(u_2(x)) = P_{W \to e}(u_2(x)). \tag{.9}$$

Then according to (.9), we get

$$\sigma_k^* = \sup_{V_k} \inf_{\|u(x)\|_{V_W}=1} \|P_{OS}^*u(x)\|_{H^{1/2}(e)} = \sup_{V_k} \inf_{\|u(x)\|_{V_W}=1} \|P_{W \to e}u(x)\|_{H^{1/2}(e)}, \tag{.10}$$

where $V_k$ is $k$-dimensional subspace of $V_W^0$, since $P_{OS}^*((V_W^0)^\perp) = 0$. And for $\sigma_k$, we have

$$\sigma_k^{\geq} \sup_{V_k} \inf_{\|u(x)\|_{V_W}=1} \|P_{OS}u(x)\|_{H^{1/2}(e)} = \sup_{V_k} \inf_{\|u(x)\|_{V_W}=1} \|P_{W \to e}u(x)\|_{H^{1/2}(e)}. \tag{.11}$$

With (.10) and (.11), we finish the proof of this Theorem.

## Proof of Lemma 3.1

Denote $f_i(x)$ as the piecewise constant functions on the coarse mesh $\Omega = \cup_{i=1}^N D_i$, then we consider $u_i(x)$ as the corresponding solution to the Poisson equation on $\Omega$. Namely,

$$-\Delta u_i(x) = f_i(x), \quad u_i(x)|_{\partial \Omega} = 0. \tag{.12}$$

Clear $u_i(x)$ are linearly independent, and we consider the projection of $u(x)$, which is the corresponding solution to the Poisson equation with $f(x)$, in the $H_0^1(\Omega)$ norm. We denote the projection as $u_h(x) = \sum_{i=1}^N c_i u_i(x)$, then we have

$$\langle u(x) - u_h(x), u_i(x) \rangle_{H_0^1(\Omega)} = 0. \tag{.13}$$

Using (.13), we have

$$\|u(x) - u_h(x)\|_{H_0^1(\Omega)}^2 = \langle u(x) - u_h(x), u(x) - u_h(x) \rangle_{H_0^1(\Omega)} = \langle u(x), u(x) - u_h(x) \rangle_{H_0^1(\Omega)}$$
$$= \langle u(x) - u_h(x), f(x) \rangle_{L^2(\Omega)} \leq \|u(x) - u_h(x)\|_{L^2(\Omega)} \|f(x)\|_{L^2(\Omega)}. \tag{.14}$$

According to (.13) and (.12), we have

$$0 = \int_\Omega (u(x) - u_h(x)) f_i(x) \mathrm{d}x = \int_{D_i} (u(x) - u_h(x)) \mathrm{d}x. \tag{.15}$$

Then based on (.15) using the Poincaré Friedrichs inequality, and a scaling argument, we have

$$\int_{D_i} (u(x) - u_h(x))^2 \mathrm{d}x \leq CH^2 \int_{D_i} |\nabla(u(x) - u_h(x))|^2 \mathrm{d}x.2\mathrm{d}x. \tag{.16}$$

Summing (.16) over $D_i$, we have

$$\|u(x) - u_h(x)\|_{L^2(\Omega)} \leq CH \|u(x) - u_h(x)\|_{H_0^1(\Omega)}. \tag{.17}$$

Putting (.17) in (.14), we have

$$\|u(x) - u_h(x)\|_{H_0^1(\Omega)}^2 \leq CH \|f(x)\|_{L^2(\Omega)} \|u(x) - u_h(x)\|_{H_0^1(\Omega)}. \tag{.18}$$

Then employing $\|\Delta(u(x) - u_h(x))\|_{H^{-1}(\Omega)} = \|u(x) - u_h(x)\|_{H_0^1(\Omega)}$, we finish the proof

$$\|f(x) - \sum_{i=1}^N c_i f_i(x)\|_{H^{-1}(\Omega)} \leq CH \|f(x)\|_{L^2(\Omega)}. \tag{.19}$$